%This is AMS-TeX

\input amstex

\magnification=1200

\documentstyle{amsppt}

\TagsOnRight

%\pagewidth{12.5cm}
%\pageheight{19.5cm}

%\hcorrection{2cm}
%\vcorrection{2cm}

\document

%Definitions (AMS-TeX)

 \def\today{\ifcase\month\or
  January\or February\or March\or April\or May\or June\or
  July\or August\or September\or October\or November\or December\fi
  \space\number\day, \number\year}

\def\om{{\omega}}
\def\sig{{\sigma}}

\def\G{{\Bbb G}}
\def\Z{{\Bbb{Z}}}
\def\Gm{{\Bbb G}_{\text{m}}}
\def\Gmk{{\Bbb G}_{{{\text{\rm{m}}},k}}}

\def\Gal{{\text{\rm{Gal}}}}
\def\Hom{{\text{\rm{Hom}}}}

\def\min{^{-1}}

%Cyrillic
\font\cyr=wncyr10 %scaled \magstep1%
\def\SH{{\text{\cyr Sh}}}

%Weil

\def\RFk{R_{F/k}}

%%%%%%%%%%%%%%%%%%%%%%%%%%%%%%%%%%%%%%%%%%%%%%%%%%%%%%%%%%%%%%%%%%%%%%%%%%%%%

\topmatter
\title On a conjecture of Le Bruyn \endtitle
\author Anne Cortella and Boris Kunyavskii \endauthor

\address A.C.: Universit\'e de Franche-Comt\'e, 
Facult\'e des Sciences et des Techniques, 
Laboratoire de Math\'ematiques, 
16 route de Gray,  
25030 BESAN\c{C}ON Cedex, 
FRANCE  \endaddress
\email cortella\@math.univ-fcomte.fr \endemail

\address B.K.: Bar-Ilan University, 
Department of Mathematics and Computer Science, 
52900 RAMAT GAN, 
ISRAEL \endaddress
\email kunyav\@macs.biu.ac.il \endemail

\date {} \enddate

\thanks 
Research of the second named author was partially 
supported by the Ministry of Absorption (Israel) and The Israel 
Science Foundation founded by the Israel Academy of Sciences and Humanities. 
\endthanks
 

\abstract Given a generic field extension $F/k$ of degree $n>3$ (i.e. the Galois group of the normal closure of $F$ is isomorphic to the symmetric group $S_n$), we prove that the norm torus, defined as the kernel of the norm map $N\colon R_{F/k}(\Bbb G_{\text{m}})\to\Bbb G_{\text{m}}$, is not rational over $k$. 
\endabstract

\endtopmatter

%%%%%%%%%%%%%%%%%%%%%%%%%%%%%%%%%%%%%%%%%%%%%%%%%%%%%%%%%%%%%%%%%%%%%%%%%%%

\document

Given an arbitrary field $k$, we call a separable 
extension $F/k$ of degree $n$ {\it generic} if the Galois group 
$G=\Gal (L/k)$ of the normal closure $L$ of $F$ over $k$ is isomorphic to 
the symmetric group $S_n$. We consider 
the norm map $N\colon F^*\to k^*$. The kernel of $N$ can be regarded as the 
set of $k$-points of an 
affine algebraic $k$-variety $T$ called {\it norm torus}. Using the Weil 
symbol of restriction of scalars, we write $T$ as the kernel of 
$\RFk (\Gm )\to\Gmk$ where $\Gm$ stands for the multiplicative group. If 
the extension $F/k$ is generic, the norm torus is also called generic and 
is denoted by $T_{F/k}$, or just $T_n$ if it does not lead to any confusion. 

In \cite{LB}, assuming $n>3$, Le Bruyn proves that the generic norm torus 
$T_n$ is non-rational over $k$ whenever $n$ is prime, and states a 
conjecture that $T_n$ is never $k$-rational except, possibly, for $n=6$. 
Our goal is to prove the above conjecture (including the case $n=6$). 
Recall that $T$ is called {\it stably rational} if there is a rational variety $T'$ such that $T\times T'$ is rational. 

\proclaim{Theorem} With the above notation, $T_n (n>3)$ is never stably  rational over $k$. 
\endproclaim 

\remark{Remark} The result might look a little bit surprising in view of good arithmetic properties of generic norm tori: in particular, if $k$ is a number field, they are known to satisfy weak approximation property and their principal homogeneous spaces satisfy the Hasse principle.  Moreover, for the case when $n$ is prime, $T_n$ is known to be a direct factor of a rational variety \cite{CT/S2}.  Note that the result cannot be ameliorated in the sense that for $n=2$ or $3$ the torus $T_n$ is of dimension 1 or 2 and hence rational \cite{V}, 4.73, 4.74.  
\endremark

\demo\nofrills{} The proof follows from the lemmas below. Throughout we denote by $M_n=\Hom (T_n,\Gm )$ the group of rational characters of $T_n$ viewed as a $G$-module. By definition, there is an exact sequence of $G$-modules 
$$
0\to\Z\to P_n\to M_n\to 0
\tag1
$$
where $P_n=\Z [G/H]$ is a permutation module, $G=S_n$, $H=\Gal (L/F)$ is isomorphic to $S_{n-1}$. The following lemma is the key one. 
\enddemo

\proclaim{Lemma 1} Let $n=rs$ with arbitrary $r,s>1$, and let $F/k$ be a generic extension of degree $n$. If $T_{F/k}=T_n$ is stably rational over $k$, there is a generic extension $K/E$ of degree $r$ such that $T_{K/E}=T_r$ is stably rational over $E$. 
\endproclaim
  
\demo{Proof} Take a subgroup $U=S_r\subset S_n$ embedded in such a way that $P_n$ restricted to $U$ is a direct sum $\underbrace{P_r\oplus\dots\oplus P_r}_{\text{$s$ times}}$. (This simply means that we partition $\{ 1,\dots ,n\}$ into $s$ disjoint subsets of cardinality $r$ and let $U$ act in a standard way on each of these subsets.) We then regard (1) as a sequence of $U$-modules and notice that $M_n$ restricted to $U$ splits into a direct sum:
$$
(M_n)|_U=M_r\oplus\underbrace{P_r\oplus\dots\oplus P_r}_{\text{$(s-1)$ times}}.  \tag2
$$
In the language of tori, (2) reads as follows: let $E=L^U$ be the fixed subfield of $U$, then the $E$-torus $T_E=T\times _kE$ is isomorphic to a direct product of $T_r=\ker [R_{K/E}(\Gm )\to\G_{\text{m},E}]$ and a quasi-split torus $S=\prod_{i=1}^{s-1}R_{K/E}(\Gm )$ where $K/E$ is a generic extension of degree $r$, $K=L^V$, $V\subset U$, $V\cong S_{r-1}$. By assumption, $T$ is stably rational over $k$, hence $T_E$ is stably rational over $E$. Since any  quasi-split torus is rational, we are done.     
\qed
\enddemo  

\proclaim{Lemma 2 (Le Bruyn)} If $n>3$ is a prime number, $T_n$ is not stably rational.  
\endproclaim

\demo{Proof} See \cite{LB}. \qed\enddemo

Before stating the next lemma, we recall that the group 
$$
\SH _{\om } ^2(G,M)=\ker [H^2(G,M)\to\prod_{g\in G}H^2(\left< g\right> ,M)]
$$
(where $M$ stands for the character module of an algebraic torus $T$ defined over $k$ and split over $L$, $G=\Gal (L/k)$), is a birational invariant of $T$. 
To be more precise, this group is zero whenever $T$ is stably rational over $k$.     Here is another useful description of the above invariant: 
consider a flasque resolution of $M$, i.e. an exact sequence of $G$-modules 
$$
0\to M\to S\to N\to 0
$$
where $S$ is a permutation module and $N$ is a flasque module (the latter means that $H\min (G',N)=0$ for all subgroups $G'\subseteq G$), then $\SH _{\om }^2(G,M)\cong H^1(G,N)$. See \cite{V}, 4.61, \cite{CT/S1} for more details. 
 
\proclaim{Lemma 3} If $n$ is a square, $T_n$ is not stably rational.
\endproclaim

\demo{Proof} Let $n=m^2$. Take a subgroup $U=\Z /m\Z\times\Z /m\Z\subset S_n$ embedded in such a way that the module $P$ from sequence (1) viewed as a $U$-module is isomorphic to $\Z [U]$. In other words, we choose $U$ generated by $$
\gather
\sig =(1\quad 2\,\dots\, m)(m+1\quad m+2\,\dots\,  2m)\dots (n-m+1\quad n-m+2\,\dots\, n), \\
\tau =(1\quad m+1 \,\dots\, n-m+1)(2\quad m+2 \,\dots\, n-m+2)\dots (m\quad 2m \,\dots\, m^2).
\endgather
$$  
Then $M_n$, regarded as a $U$-module, is none other than $\hat J=\Z [U]/\Z$, the character module of the norm torus $J=\ker [R_{L/E}(\Gm )\to\G_{\text{m},E}]$ where $E=L^U$. It is well known that $J$ is not rational over $E$ because 
$\SH _{\om }^2(U,\hat J)=\Z /p\Z$. Since $T_n\times _kE=J$, we conclude that 
$T_n$ cannot be stably rational over $k$. \qed
\enddemo

\proclaim{Corollary (Saltman, Snider)} If $n$ is divisible by a square, $T_n$ is not stably rational.
\endproclaim

\demo{Proof} Combine Lemma 1 and Lemma 3. \qed
\enddemo

\proclaim{Lemma 4} The torus $T_6$ is not stably rational. 
\endproclaim

\demo{Proof} Take $U=\Z /2\times\Z /2\subset S_6$ generated by $(12)(34)$ and 
$(34)(56)$. We observe that $U$ coincides with the Sylow 2-subgroup of the 
alternating group $A_4$ embedded into $S_6$ via its action on the edges of 
tetrahedron.  Let $M=M_6$ be the module of characters of $T_6$ defined by 
sequence (1) with $G=S_6$, $H=S_5$. It is known that 
$\SH _{\om }^2(A_4,M)=\Z /2\Z$ (\cite{D/P}, Lemma 13). This implies 
$\SH _{\om }^2(U,M)\neq 0$. Indeed, assume the contrary. Then, since any 
Sylow 3-subgroup $V$ of $A_4$ is cyclic, one has $\SH _{\om }^2(V,M)=0$, 
and vanishing of $\SH _{\om }^2(U,M)$ would imply vanishing of 
$\SH _{\om }^2(A_4,M)$ (one may apply the above interpretation of 
$\SH _{\om }^2(G,M)$ as $H^1(G,N)$ to the case $G=A_6$ and use the fact 
that the restriction to a Sylow $p$-subgroup is injective on the 
$p$-component of $H^1$). 
\qed
\enddemo

\remark{Remark} Of course, one may give a more direct proof of Lemma 4 without referring to \cite{D/P}, either by a straightforward computation of $\SH _{\om }^2(U,M)$ (which goes much simpler than for $A_4$), or by constructing an exact sequence of $U$-modules 
$$
0\to M_a\to M\to \Z\oplus\Z\to 0
$$
with $M_a$ the character module of an anisotropic torus $T_a$ which, in our case, turns out to be $\Z[U]/\Z$; by \cite{V}, 4.22, the latter exact sequence induces a birational equivalence of tori $T_n\sim T_a\times\G ^2_{\text{m}}$, whence the result. 
\endremark

\demo{Proof of the Theorem} We are now ready to prove the Theorem. Indeed, the above Corollary reduces the problem to the case when $n$ is square-free, and Lemmas 1 and 2 englobe all $n$ having a prime divisor greater than 3. We thus have to apply Lemma 4 for the only remaining case $n=6$. 
\qed
\enddemo

\remark{Concluding remark} Our theorem can (and should) be viewed in a broader context. Namely, one can extend it to generic tori in (almost absolutely) simple groups. Indeed, the above result corresponds to the case of an inner form of a simply connected group of type $A_{n-1}$. Such a generalization to the other types of inner and outer forms of simply connected and adjoint groups is the subject of our forthcoming paper. 
\endremark

\remark{Acknowledgement} This work was done while the second named author was visiting Universit\'e de Franche-Comt\'e at Besan\c{c}on. B. Kunyavskii thanks this institution for its hospitality. 
\endremark

%%%%%%%%%%%%%%%%%%%%%%%%%%%%%%%%%%%%%%%%%%%%%%%%%%%%%%%%%%%%%
\refstyle{A}

%%%%%%%%%%%%%%%%%%%%%%%%%%%%%%%%%%%%%%%%%%%%%%%%%%%%%%%%%%%%%%%%%%%%%%%%%%%%%%%

\enddocument